\magnification=\magstep1
\input amstex
\documentstyle{amsppt}
\hsize=6.54truein
 
\ifnum\pageno=1\nopagenumbers\fi
\nologo

\def\ges{{\geqslant}}
\def\ti{{t^{-1}}}
\def\QED{\hfill\hbox{\qed}}

\def\cx{^{\ssize\bullet}}

\def\n{{\frak n}}
\def\p{{\frak p}}
\def\q{{\frak q}}
\def\r{{\frak r}}
\def\s{{\frak s}}
\def\f{{\varepsilon}}
\def\BA{{\Bbb A}}
\def\BC{{\Bbb C}}
\def\BZ{{\Bbb Z}}

\def\E#1{{}_{#1}\!\operatorname{E}}
\def\e#1{{}^{#1}\!\operatorname{E}}
\def\EE#1{{}_{#1}\!\operatorname{E}^{q}_{p}}
\def\Ee#1{{}^{#1}\!\operatorname{E}^{p,q}}

\def\Supp{\operatorname{Supp}}
\def\Tor{\operatorname{Tor}}
\def\Ext{\operatorname{Ext}}
\def\Hom{\operatorname{Hom}}
\def\Ker{\operatorname{Ker}}
\def\Im{\operatorname{Im}}
\def\Proj{\operatorname{Proj}}
\def\rank{\operatorname{rank}}
\def\deg{\operatorname{deg}}
\def\ord{\operatorname{ord}}
\def\depth{\operatorname{depth}}
\def\height{\operatorname{height}}
\def\codim{\operatorname{codim}}
\def\length{\operatorname{length}}
\def\pd{\operatorname{proj\,dim}}
\def\id{\operatorname{inj\,dim}}
\def\G{\operatorname{G}}
\def\ann{\operatorname{ann}}
\def\lau{{f}_R}

\topmatter
\title
Laurent coefficients and Ext of finite graded modules
\endtitle

\leftheadtext{L.~L.~Avramov, R.-O.~Buchweitz, and J.~D.~Sally}
\rightheadtext{Laurent coefficients and Ext}

\author
Luchezar~L.~Avramov\\
Ragnar-Olaf~Buchweitz\\
Judith~D.~Sally
\endauthor

\address
Department of Mathematics, Purdue University, West Lafayette,
Indiana 47907, U.~S.~A.
\endaddress
\email
{avramov\@math.purdue.edu}
\endemail

\address
Department of Mathematics, University of Toronto,
Toronto, Ontario M5S 1A1, Canada
\endaddress
\email
{ragnar\@lake.scar.utoronto.ca}
\endemail

\address
Department of Mathematics, Northwestern University,
Evanston, Illinois 60208, U.~S.~A.
\endaddress
\email
{judy\@math.nwu.edu}
\endemail

\endtopmatter

\subjclass
13D40, 13F15, 13H10
\endsubjclass

\thanks
Research partially supported by NSF (L.~L.~A. and J.~D.~S.)
and NSERC (R.-O.~B.).
\endthanks
\endtopmatter

\document

Let $R=\bigoplus_{n\ges0}R_n$ be a graded commutative ring
generated over a field $K=R_0$ by homogeneous elements $x_1,\dots,x_e$
of positive degrees $d_1,\dots,d_e$. The Hilbert-Serre Theorem shows that
for each finite graded $R$--module $M=\bigoplus_{n\in\BZ}M_n$ the {\it Hilbert
series\/} $\sum_{n\in\BZ}(\rank_K M_n)t^n$ is the Laurent expansion around
$0$ of a rational function
$$
H_M(t)=\frac{q_M(t)}{\prod_{i=1}^e(1-t^{d_i})}
$$
with $q_M(t)\in\BZ[t,\ti]$.  We demonstrate that Laurent expansions
$\left[M\right]_z$ of $H_M(t)$ around other points $z$ of the extended complex
plane $\overline\BC$ also carry important structural information.

When $R$ is generated in degree one and $z=1$ such an approach has been used
implicitly, since the coefficients of the principal part of $\left[M\right]_1$
are then also coefficients of the Hilbert polynomial of $M$.  Without
assumptions on the generation of $R$, coefficients of $\left[M\right]_1$ have
been sighted occasionally in invariant theory.  In that context Benson and
Crawley-Boevey \cite{3} have discovered recently that for all finite
graded modules over a graded normal domain the coefficient $\psi$ of
$1/(1-t)^{\dim R-1}$ satisfies an equality
$$
\gathered
\psi(\Hom_R(M,N))-\psi(\Ext^1_R(M,N))=\\
\rank_R(M)\rank_R(N)\psi(R)
+\rank_R(M)\psi(N)-\rank_R(N)\psi(M)\,.
\endgathered\tag1
$$

The present investigation started as an attempt to ``explain'' this intriguing
formula, and to find out whether analogous equalities exist for subsequent
coefficients of Laurent expansions around $1$.  Fairly complete answers are
presented in Sections 3 and 4.  They depend on embedding (1) into an equality
of Laurent series, which ``predicts'' an infinite sequence of similar
relations among higher Laurent coefficients.  The results are proved in
Section 1, and are used in Section 2 to obtain non-trivial lower bounds for
the Bass numbers of certain graded modules of finite injective dimension. 

\head
1. Hilbert series
\endhead

When $M$ and $N$ are graded $R$--modules $\Hom_R(M,N)_n$ stands for the
$K$--vector subspace of $\Hom_R(M,N)$ consisting of homomorphisms $\alpha$ of
{\it degree\/} $n$, that is, such that $\alpha(M_j)\subseteq N_{n+j}$ for
$j\in\BZ$.  For a finite $M$ the inclusion $\bigoplus_{n\in\BZ}\Hom_R(M,N)_n
\subseteq\Hom_R(M,N)$ is an equality, and thus $\Hom_R(M,N)$ is graded; as $M$
has a resolution by finite free graded $R$--modules and homomorphisms of degree
$0$, this induces for each $i$ a grading on $\Ext^i_R(M,N)$.

\proclaim{Theorem 1}
If $M$ and $N$ are finite graded $R$--modules such that $\Ext^i_R(M,N)=0$
for $i\gg0$, then there is an equality of rational functions
$$
\sum_{i}(-1)^iH_{\Ext^i_R(M,N)}(t)=\frac{H_M(\ti)\cdot H_{N}(t)}{H_R(\ti)} \,.
$$
\endproclaim

\remark{Remark}
This is a far reaching generalization of the known expression 
$H_{\Ext^{d-m}_R(M,\omega)}(t)=(-1)^mH_M(\ti)$
which holds for an $m$-dimensional graded Cohen-Mac\-aulay module $M$ over a
$d$-dimensional Cohen-Macau\-lay ring $R$ with canonical module $\omega$,
cf\.  \cite{4\rm; (4.3.7)}.

Indeed, $\Ext^i_R(M,\omega)=0$ unless $i=d-m$ by \cite{4\rm; (3.3.10)}, so
we only need to show that $H_{\omega}(t)=(-1)^dH_R(\ti)$.  Consider the
polynomial ring $Q=K[X_1,\dots,X_e]$ graded by $\deg X_j=d_j$ and the
surjective homomorphism of $K$--algebras $Q@>>>R$ with $X_j\mapsto x_j$ for
$1\le j\le e$.  By \cite{4\rm; (3.6.12), (3.6.10)}, for $W=Q(-\sum_{i=1}^ed_i)$
we have $\omega\cong\Ext^{e-d}_Q(R,W)$, hence we get $H_{\omega}(t)$ by
applying the theorem to the $Q$--modules $R$ and $W$.
\endremark
\medskip

When $f(t)$ is a function of a complex variable $t$, we denote by
$\left[f(t)\right]_z$ the expansion of $f(t)$ as a Laurent series
$\sum_{j\in\BZ}a_j(t-z)^j$ when $z\in\BC$, and as a Laurent series
$\sum_{j\in\BZ}a_j t^{-j}$ when $z=\infty$. The {\it order\/} of
$\left[f(t)\right]_z$ is introduced by
$\ord\left[f(t)\right]_z=\inf\{j\in\BZ\mid a_j\ne0\}$.

For a finite graded $R$--module $M$ we often write $\left[M\right]_z$ in place
of $\left[H_M(t)\right]_z$.  Clearly, $\left[M\right]_0\in\BZ((t))$ and
$\left[M\right]_{\infty}\in\BZ((\ti))$, where $\BA((t))=\BA[[t]][\ti]$ and
$\BA((\ti))=\BA[[\ti]][t]$ denote rings of formal Laurent series of finite
order with coefficients in $\BA$, with topologies defined by the powers of
the respective indeterminate.  Note that if $\{g_i\}_{i\ges0}$ is a sequence of
formal Laurent series in either ring, then $\sum_{i\ges0}g_i$ converges when
$\lim_{i\to\infty}\ord g_i=\infty$.

Much of the paper deals with properties of the rational function
$$
\phi_R(M,N)(t)=\frac{H_M(\ti)\cdot H_{N}(t)}{H_R(\ti)}\,.
$$
Without homological assumptions on $M$ or $N$, we establish the following
limited

\proclaim{Proposition 2}
If $M$ and $N$ are finite graded $R$--modules, and $N$ has finite length,
then $H_{\Ext^i_R(M,N)}(t)\in\BZ[t,\ti]$,
the order of $\left[\Ext^i_R(M,N)\right]_\infty$ goes to infinity
with $i$, and $\sum_i(-1)^i\left[\Ext^i_R(M,N)\right]_\infty\in\BZ((\ti))$
is the Laurent expansion around $\infty$ of $\phi_R(M,N)(t)$.
\endproclaim

For comparison and for later application, we recall the situation in homology.
The grading of $M\otimes_KN$ by
$(M\otimes_KN)_n=\bigoplus_{j\in\BZ}(M_j\otimes_KN_{n-j})$ induces a grading
of $M\otimes_RN$. As above, this defines a grading of $\Tor^R_i(M,N)$ for each
$i$, and by  \cite{1\rm; Lemma 7} we have

\proclaim{Lemma 3}
If $M$ and $N$ are finite graded $R$--modules, then the order of
$\left[\Tor^R_i(M,N)\right]_0\in\BZ((t))$ goes to $\infty$ with $i$, and
$\sum_{i}(-1)^i\left[{\Tor^R_i(M,N)}\right]_0\in\BZ((t))$
is the Laurent expansion around $0$ of the rational function
\TagsOnRight
$$
\chi^R(M,N)(t)=\frac{H_M(t)\cdot H_N(t)}{H_R(t)}\,.\tag"\qed"
$$
\TagsOnLeft
\endproclaim

Before starting on the proofs we record some general properties of
Laurent expansions.

\remark{Remark}
If $0@>>>M'@>>>M@>>>M''@>>>0$ is an exact sequence of graded $R$--modules, then
$\left[M'\right]_z+\left[M''\right]_z=\left[M\right]_z$.  Indeed, additivity
is obvious for $z=0$, hence yields an equality of rational functions
$H_{M'}(t)+H_{M''}(t)=H_{M}(t)$, which implies additivity for arbitrary $z$.

Similarly, one sees that
$\left[M(a)\right]_z=\left[t^{-a}\right]_z\left[M\right]_z$ for any $a\in\BZ$,
where $M(a)$ is the $a$'th {\it translate\/} of $M$, that is, the graded
$R$--module with $M(a)_n=M_{a+n}$ for $n\in\BZ$.
\endremark

\demo {Proof of Theorem \rom{1}}
Let $Q@>>>R$ and $W$ be as in the remark following the statement of the theorem,
let $I\cx$ be a minimal (hence finite, cf\. \cite{4\rm, (3.6.6)}) graded
injective resolution of the $Q$--module $W$, and let $F\cx$ be a minimal graded
free resolution of the $R$--module $M$. As in \cite{2\rm, (4.4.I)}, there is a
canonical isomorphism of complexes of graded $R$--modules
$$
\theta\:F\cx\otimes_R\Hom_Q(N,I\cx)@>>>\Hom_Q(\Hom_R(F\cx,N),I\cx)
$$
given by $\theta(f\otimes_R\alpha)(\beta)=(-1)^{i(j-h)}\alpha\beta(f)$ for
$f\in F^{i}$, $\alpha\in\Hom_Q(N,I^j)$, $\beta\in\Hom_Q(F^{h},N)$.
Filtering the complex on the left by the cohomological degree of $F\cx$, and
the one on the right by that of $I\cx$, we obtain two spectral sequences
which converge to a common limit:
$$
\EE{2}=\Tor^R_{p}(M,\Ext^{q}_Q(N,W))\Longrightarrow E_{p-q}
\Longleftarrow\Ext^{q}_Q(\Ext^p_R(M,N),W)=\Ee{2}\,.
$$

In the first spectral sequence we have
${}_rd^q_p\:\EE{r}@>>>\E{r}^{q-r+1}_{p-r}$ for $r\ge2$, and $\EE{2}=0$ unless
$p\ge0$ and $0\le q\le e$.  It follows that $E_i=0$ for $i<-q$, and that
$\EE{e+2}=\EE{\infty}$.

By Lemma 3 we know that for any fixed $q$ the order of
$\left[{\Tor^R_{p}(M,\Ext^{q}_Q(N,W))}\right]_0\in\BZ((t))$ goes to
infinity together with $p$, hence $\sum_p(-1)^p\left[\EE{2}\right]_0$
is in $\BZ((t))$. As $\EE{r}$ is a subquotient of $\EE{2}$, we further have
$\sum_p(-1)^p \left[\EE{r}\right]_0\in\BZ((t))$ for each $r\ge2$, and
we set
$$
\chi(\E{r})=\sum_{q=0}^e(-1)^q\sum_p(-1)^p\left[{\EE{r}}\right]_0\,.
$$

Next we show that $\chi(\E{r})=\chi(\E{r+1})$ for $r\ge2$. To this
end, we first note that up to any given degree $n$ each series
is determined by a finite number of the summands used to define it, hence when
looking at the coefficient of $t^n$ we may restrict both summations to the
same finite set of indices $p$ and $q$.  In $\chi(\E{r})$ this coefficient
is the Euler characteristic of the degree $n$ part of a finite complex of
graded $R$--modules.  As the differential ${}_rd^q_p$ preserves the grading
of these modules, the coefficient of $t^n$ in $\chi(\E{r+1})$ is the Euler
characteristic of the degree $n$ part of the homology of this complex.  The
classical permanence property of Euler characteristics implies that these
coefficients are equal.

By the finite convergence of the spectral sequence, we thus see that
the series $\chi(E)=\sum_i(-1)^i\left[{E_i}\right]_0$
is equal to $\chi(\E{2})$. Together with Lemma 3 this implies
$$
\align
\chi(E)
&=\sum_{q=0}^e(-1)^q
\sum_p(-1)^p\left[{\Tor^R_p(M,\Ext^q_Q(N,W))}\right]_0\\
&=\sum_{q=0}^e(-1)^q\left[\chi_R(M,\Ext^q_Q(N,W))(t)\right]_0\\
&=\left[\frac{H_M(t)}{H_R(t)}\sum_{q=0}^e(-1)^qH_{\Ext^q_Q(N,W))}(t)\right]_0\,.
\endalign
$$

In order to finish the computation we use the equality of rational functions
$$
\sum_{q=0}^{e}(-1)^q H_{\Ext^q_Q(A,W)}(t)=(-1)^eH_A(\ti)\,.\tag2
$$
which holds for each finite $Q$--module $A$.  Indeed, if $A=Q(b)$, then the
only non-vanishing $\Ext$ is $\Ext^0_Q(Q(b),W))\cong Q(-b-\sum_{j=1}^ed_j)$,
and the equality is checked by a direct computation which uses the expression
$H_Q(t)=1/\prod_{j=1}^e(1-t^{d_j})$.  The general case follows, as each
$A$ has a finite resolution by finite direct sums of translates
of $Q$, and both sides of the formula are additive functions on the
Grothendieck group $\G(Q)$ of the category of finite graded $Q$--modules and
homomorphisms of degree $0$.  Thus, we have established that
$$
\chi(E)=(-1)^e\left[\frac{H_M(t)\cdot H_N(\ti)}{H_R(t)}\right]_0
=(-1)^e\left[\phi_R(M,N)(\ti)\right]_0\,.
$$

Now we turn to the second spectral sequence. It has $\Ee{2}=0$ unless $p\ge0$
and $0\le q\le e$, and differentials ${}^rd^{p,q}\:\Ee{r}@>>>\e{r}^{p-r+1,q-r}$,
hence $\Ee{e+2}=\Ee{\infty}$.  Furthermore, our assumption implies that
$\Ee{r}=0$ for $p\gg0$, hence for each $r\ge2$ the sum 
$$
\chi(\e{r})=\sum_p(-1)^p\sum_{q=0}^e(-1)^q\left[{\Ee{r}}\right]_0
$$
is actually {\it finite\/}.  This time we may apply directly the
classical argument on Euler characteristics.  In view of formula (2) it yields
$$
\align
\chi(E)
&=\sum_p(-1)^p\sum_{q=0}^e(-1)^q\left[{\Ext^q_Q(\Ext^p_R(M,N),W)}\right]_0\\
&=\sum_p(-1)^p
\left[\sum_{q=0}^e(-1)^qH_{\Ext^q_Q(\Ext^p_R(M,N),W)}(t)\right]_0\\
&=(-1)^e\left[\sum_p(-1)^pH_{\Ext^p_R(M,N)}(\ti)\right]_0\,.
\endalign
$$

The overall result of the preceding computations now reads
$$
\left[\shave\sum_p(-1)^pH_{\Ext^p_R(M,N)}(\ti)\right]_0=
\left[\phi_R(M,N)(\ti)\right]_0
$$
and this clearly implies the desired equality of rational functions.
\QED
\enddemo

\demo{Proof of Proposition \rom{2}}
Consider $N^\vee=\Hom_K(N,K)$ with the induced structure of graded $R$--module.
By Lemma 3, the expansion around $0$ of the rational function $\chi^R(M,N)(t)$
is equal to $\sum_i(-1)^i\left[\Tor^R_i(M,N^\vee)\right]_0\in\BZ((t))$.
As $\rank_KN$ is finite, we have $H_{N^\vee}(t)=H_N(\ti)\in
\BZ[t,\ti]$, and thus an equality of formal Laurent series
$$
\sum_{i\ges0}(-1)^i\left[H_{\Tor^R_i(M,N^\vee)}(t)\right]_0=
\left[\frac{H_M(t)\cdot H_{N}(\ti)}{H_R(t)}\right]_0\,.
$$

In the commutative diagram of homomorphisms of rings
$$
\CD
\BC(t)@>t\mapsto\ti>>\BC(t)\\
@V\left[\ \right]_0 VV @VV\left[\ \right]_\infty V\\
\BC((t))@>>t\mapsto\ti>\BC((\ti))
\endCD
$$
the lower row is an isomorphism of topological fields, hence the preceding
equality yields
$$
\sum_{i\ges0}(-1)^i\left[H_{\Tor^R_i(M,N^\vee)}(\ti)\right]_\infty=
\left[\frac{H_M(\ti)\cdot H_{N}(t)}{H_R(\ti)}\right]_\infty\,.
$$
The canonical isomorphisms $(M\otimes_R(N^\vee))^\vee\cong\Hom_R(M,N^{\vee\vee})
\cong \Hom_R(M,N)$ extend to isomorphisms of graded $R$--modules
$(\Tor^R_i(M,N^\vee))^\vee\cong\Ext^i_R(M,N)$ for $i\ge0$. Thus, in
$\BZ[t,\ti]$ we have $H_{\Tor^R_i(M,N^\vee)}(\ti)=H_{\Ext^i_R(M,N)}(t)$,
and this finishes the proof.
\QED
\enddemo

\head
2. Bass numbers
\endhead

Let $N$ be a finite graded $R$--module.  For an integer $i$ and a prime ideal
$\p$ of $R$, the rank of $\Ext^i_{R_\p}(k(\p),N_\p)$ over the field
$k(\p)=R_\p/\p R_\p$ is known as the $i$'th {\it Bass number\/} of $N$ and is
denoted $\mu^i_R(\p,N)$.  We establish nontrivial lower
bounds for the Bass numbers $\mu^i_R(N)=\mu^i_R(R_+,N)$ of ``many'' graded
modules of finite injective dimension.

A graded ring $R$ is said to be {\it standard\/} if it is generated by elements
of degree $1$.  In this case $H_N(t)$ has a unique expression of the form
$$
H_N(t)=\frac{e_N(t)}{(1-t)^n}\,,\tag3
$$
where $n=\dim N$, and $e_N(t)$ is a Laurent polynomial in $\BZ[t,\ti]$ such
that $e_N(1)$ is a positive integer, known as the {\it multiplicity\/} of $N$.  

\proclaim{Theorem 4}
Let $R$ be a standard graded ring of Krull dimension $d$.

If $N\ne0$ is a finite graded $R$--module of finite injective dimension and 
Krull dimension $n$, then $e_R(\ti)$ divides $e_N(t)$ in $\BZ[t,\ti]$.
If $p$ is a prime number and $q$ is defined by
$$
q=\max\,\biggl\{p^r\in\BZ\biggl.\ \overset\mathstrut\to{\mathstrut}\biggr|_{\,}
\,\sum_{s=0}^{p^r}t^s \text{\ \,divides\ \,} \frac{e_N(t)}{e_R(\ti)}
\text{\ \,in\ \,} \BZ[t,\ti]\biggr\}\,,
$$
then the Bass numbers of $N$ satisfy the inequality
$$
\sum_i\mu^i_R(N)\ge p^{\tsize{\frac{d-n+q-1}{q(p-1)}}}\,.
$$
\endproclaim

\demo{Proof}
Since each $\Ext^i_R(K,N)$ is a finite graded $K$--vector space, and is trivial
for $i>\id_RN$, the rational function $\sum_i(-1)^iH_{\Ext^i_R(K,N)}(t)$ is a
Laurent polynomial, which we denote by $\phi^N_R(t)$.  Theorem 1 and formula
(3) yield an equality in $\BZ[t,\ti]$:
$$
(1-t)^{d-n}e_N(t)=(-t)^de_R(\ti)\phi^N_R(t)\,.
$$
As $d\ge n$ and $e_R(1)\ne0$, unique factorization in $\BZ[t,\ti]$ shows
that $e_R(\ti)$ divides $e_N(t)$.  This proves the first assertion.  The
second one follows by the argument of \cite{1\rm; \S4}.
\QED
\enddemo

Under the hypotheses of the next result the theorem applies with $p=2$ and
$q=1$:

\proclaim{Corollary 5}
If $-1$ is not a root of the Laurent polynomial $e_N(t)/e_R(\ti)$, then
$$
\sum_i\mu^i_R(N)\ge2^{d-n}\,.
$$ 
This inequality holds, in particular, when the multiplicity of $N$ is an odd
number.
\QED
\endproclaim

\remark{Remark}
If $N$ is a finite module over a local ring $(S,\n)$, and
$\depth_SN=g<d=\dim S$, then a result of Bruns \cite{4\rm; (9.6.1.a)}
implies that $\sum_i\mu^i_S(N)\ge (d-g)(d-g-1)+3$.

The {\it quadratic\/} lower bound above is in general
weaker than the {\it exponential\/} bound given by the corollary, so we raise
the question whether an inequality $\sum_i\mu^i_S(\n,N)\ge2^{d-n}$ holds for
each finite $S$--module $N$ of Krull dimension $n$ and of {\it finite injective dimension\/}.  Such a bound would be best possible: if $\bold x$ is a system
of parameters of $S$, $E$ is the injective envelope of $S/\n$, and $N$ is the
finite length module $\Hom_S(S/(\bold x),E)$, then $\mu^i_S(\n,N)=\binom di$
since $S$ is Cohen--Macaulay by the Bass Conjecture proved by Peskine, Szpiro,
Hochster, and P.~Roberts, cf\. \cite{4\rm; (9.6.2), (9.6.4.b)}.
\endremark

\head
3. Laurent coefficients
\endhead

We turn to properties of the coefficients of Laurent expansions around $1$ of
rational functions representing the Hilbert series of finite $R$--modules.  For
purposes of calibration, it is convenient to write such an expansion in the form
$$
\left[M\right]_1=\sum_{j\ges0}\frac{\lau^{j}(M)}{(1-t)^{d-j}}
$$
with $d$ equal to the Krull dimension of $R$, and to set $\lau^{j}(M)=0$ for
$j<0$.  We call $\lau^j(M)$ the $j$'th {\it Laurent coefficient\/} of $M$ and
note that it is a rational number.

\remark{Remark}
When $R$ is standard and $d-\dim M=h$, comparison with (3) yields
$\rank_K(M_n)=\sum_{j=h}^{d-1}\lau^j(M)\binom{n+d-j-1}{d-j-1}$ for $n\gg0$,
hence the Laurent coefficients of the principal part of $\left[M\right]_1$ are
directly related to the coefficients of the Hilbert polynomial of $M$ in the
usual ``binomial'' format, as in \cite{4\rm; (4.1.5)}:
$\lau^j(M)=(-1)^{h-j}e_{h-j}(M)$ for $h\le j\le d-1$.
\endremark
\medskip

Under the assumptions of Theorem 1, the $j$'th coefficient in the Laurent
expansion around $1$ of the rational function $\phi_R(M,N)(t)$ is an
alternating sum of the corresponding Laurent coefficients of $\Ext$ modules. 
The number
$$
\f_R^j(M,N)=\sum_{i=0}^j(-1)^{i}\lau^j(\Ext_R^i(M,N))
$$
represents a portion of such a sum.  The next result shows that under a
homological finiteness hypothesis it contains the entire information.

\proclaim{Proposition 6}
If the projective dimension of a finite graded $R$--module $M$ or the
injective dimension of a finite graded $R$--module $N$ is finite,
then
$$
\left[\phi_R(M,N)(t)\right]_1= \sum_{j\ges0}\frac{\f_R^j(M,N)}{(1-t)^{d-j}}\,.
$$
\endproclaim

The proposition applies to all finite modules over a graded polynomial
ring $R$, and raises the question of extending its validity---at least on
part---to other graded rings.  When
$$
\phi_R(M,N)(t)=\sum_{j=0}^n\frac{\f_R^j(M,N)}{(1-t)^{d-j}}+
O\left(\frac1{(1-t)^{d-n-1}}\right)
$$
for some integer $n\ge0$, we say that the sequences $\phi^*_R(M,N)$ (of
coefficients of the Laurent expansion of $\phi_R(M,N)(t)$ around $1$) and 
$\f_R^*(M,N)=\{\f_R^j(M,N)\}_{j\in\BZ}$ {\it agree up to level\/} $n$.

The next result shows that an agreement of $\f^*$ and $\phi^*$ imposes strong
restrictions on the local structure of $R$.  We adopt the convention that
the {\it codimension\/} of a prime ideal $\p$ of $R$ is the non-negative integer
$\codim\p=\dim R -\dim(R/\p)$ and note that $\height\p\le\codim\p$, with
equality when $R$ is an integral domain.  The ring $R$ is {\it regular in
codimension\/} $c$ if $R_\p$ is regular for each prime $\p$ such that
$\codim\p\le c\,$; by Matijevic \cite{8\rm; (2.1)} it suffices to
impose the condition only on the primes in $\Proj R$, the set of
{\it homogeneous\/} prime ideals of $R$.

\proclaim{Proposition 7}
Let $\p$ be a homogeneous prime ideal of $R$ with $\codim\p=h$.  The sequences
$\f_R^*(R/\p,R/\p)$ and $\phi^*_R(R/\p,R/\p)$ agree up to level $h$ precisely
when the local ring $R_\p$ is regular, and if furthermore $h=0$ then $\p$ is
the unique codimension $0$ prime of $R$.
\endproclaim 

The condition that $\f_R^*(M,N)$ and $\phi^*_R(M,N)$ agree up to level $n$
can be rewritten as
$$
H_R(\ti)\cdot\sum_{j=0}^n\frac{\f_R^j(M,N)}{(1-t)^{d-j}}
=H_M(\ti)\cdot H_N(t)+O\left(\frac1{(1-t)^{2d-n-1}}\right)\,.
$$
Taking on both sides Laurent expansions around $1$, we can further reformulate
it in terms of a system of $n+1$ numerical equalities involving the numbers
$\lau^j(M)$, $\lau^j(N)$, $\lau^j(R)$, and $\f_R^j(M,N)$ for $0\le j\le n$.
Equations (4.0), (4.1), (4.2) display such a system for $n=2$.

\proclaim{Theorem 8}
Let $M$ and $N$ be finite graded modules over a graded ring $R$.

If $R$ has a unique prime $\r$ of codimension $0$ and $R_\r$ is a field, then
$$
\lau^0(R)\f_R^0(M,N)=\lau^0(M)\lau^0(N)\,.\tag4.0
$$

If $R$ is an integral domain which is regular in codimension $1$, then
$$
\lau^0(R)\f_R^1(M,N)-\lau^1(R)\f_R^0(M,N)=
\lau^0(M)\lau^1(N)-\lau^1(M)\lau^0(N)\,.\tag4.1
$$

If $R$ is a unique factorization domain which is regular in codimension $2$,
then
$$
\gathered
\lau^0(R)\f_R^2(M,N)-\lau^1(R)\f_R^1(M,N)+
\left(\lau^2(R)-\lau^1(R)\right)\f_R^0(M,N)=\\
\lau^0(M)\lau^2(N)-\lau^1(M)\lau^1(N)+
\left(\lau^2(M)-\lau^1(M)\right)\lau^0(N)\,.
\endgathered\tag4.2
$$
\endproclaim

\remark{Remark}
When $R$ is a domain for any finite graded $R$--module $M$ one has
$$
\lau^0(M)=\rank_R(M)\lau^0(R)\,,\tag5
$$
cf\. Lemma 10(c) below.  If furthermore $R$ is regular in codimension $1$, then
formula (1), due to Benson and Crawley-Boevey \cite{3\rm; (2.4)}, is
seen to follow from (4.0) and (4.1).
\endremark
\medskip

Recall that a {\it graded complete intersection\/} is a quotient of a graded
polynomial ring by a regular sequence of homogeneous elements.  By a classical
result of Grothendieck \cite{7\rm; (XI.3.14)}, cf\. also \cite{5}, a complete
intersection which is regular in codimension $3$ is factorial, hence we get the
following easily applicable

\proclaim{Corollary 9}
If $R$ is a graded complete intersection which is regular in codimension
$3$, then the equalities of the theorem hold for all finite graded $R$--modules.
\QED
\endproclaim

In the next section we show that neither the UFD hypothesis of the theorem nor
the codimension $3$ hypothesis of the corollary can be significantly weakened.
The proofs of the results above depend on a few lemmas; the first one
(for $\psi=\lau^1$) is in \cite{3\rm; (2.2), (2.3)}.

\proclaim{Lemma 10}
Let $M\ne0$ be a finite $R$--modules with $h=\dim R-\dim M$.

\rom{(a)} $\lau^j(-)$ is an additive function on $\G(R)$ for each $j\in\BZ$.

\rom{(b)} $\lau^j(M)=0$ for $j<h$ and $\lau^h(M)>0$.

\rom{(c)} $\lau^{h}(M)=\sum\Sb \p\in\Proj R\\ \codim\p=h\endSb
\length_{R_\p}\!(M_\p)\lau^{h}(R/\p)$\,.
\endproclaim

\demo{Proof}
(a) results from the additivity of $\left[\ \right]_1$ and the uniqueness of
Laurent expansions.

(b) is a fact of dimension theory.

(c) is a formal consequence of the preceding two because
$M$ admits a finite filtration with subquotients of the form $(R/\q)(a)$
for appropriate $\q\in\Proj R$ and $a\in\BZ$.
\QED
\enddemo

\proclaim{Lemma 11}
Let $M$ and $N$ be finite $R$--modules, such that for some $c\in\Bbb N$ and
for all $\p\in\Proj R$ with $\codim\p\le c$ either $\pd_{R_\p}(M_\p)$ or
$\id_{R_\p}(N_\p)$ is finite.  For each integer $j\le c$ and any
submodule $L$ of $\Ext_R^i(M,N)$, we have $\lau^j(L)=0$ when $i>j$.
\endproclaim

\demo{Proof}
By Lemma 10(b) it suffices to show that $j<\dim R-\dim L=g$, hence we may assume
that $g\le c$.  If $\p\in\Supp L$ is a homogeneous prime with $\codim\p=g$, then
$$
\Ext^{i}_{R_\p}(M_\p,N_\p)\cong\Ext^{i}_R(M,N)_\p\supseteq L_\p\ne0
$$
implies that $i\le\min\{\pd_{R_\p}(M_\p),\,\id_{R_\p}(N_\p)\}$. One of these
dimensions is finite by assumption, and thus is bounded above
by $\depth R_\p$, cf\. \cite{4\rm; (1.3.3) and (3.1.17)}.
Thus, we have $j<i\le \depth R_\p\le\codim\p=g$, as desired.
\QED
\enddemo

\demo\nofrills{Proof of Proposition \rom{6}}\,:\ 
Equate the Laurent expansions around $1$ of the rational functions in
Theorem 1 and use Lemma 11 to regroup the terms in the sum of $\Ext$'s.
\QED
\enddemo

The argument for part (a) of the next lemma is abstracted from the proof of
\cite{3\rm; (2.4)}.

\proclaim{Lemma 12}
The numbers $\f_R^j(M,N)$ have the following properties.

\rom{(a)}
$\f_R^j(-,-)$ is a biadditive function on $\G(R)$ if $R$ is regular in
codimension $c$ and $j\le c$.

\rom{(b)}
$\f_R^j(M,N)=0$ when $j<\dim R-\dim(R/(\ann M +\ann N))$\,.
\endproclaim

\demo{Proof}
(a) Fix an integer $j$ between $0$ and $c$, take a short exact sequence of
finite $R$--modules $0@>>>N'@>>>N@>>>N''@>>>0$, and consider the exact sequence
$$
\alignedat7
0&@>>>&&\Hom_R(M,N')&&@>>>&&\Hom_R(M,N)&&@>>>&&\Hom_R(M,N'')&&@>>>\\
&&&\Ext^1_R(M,N')&&@>>>&&\quad\qquad\dots&&@>>>&&\Ext^{j-1}_R(M,N'')&&@>>>\\
&&&\Ext^j_R(M,N')&&@>>>&&\Ext^j_R(M,N)&&@>>>&&\Ext^j_R(M,N'')&&
@>>>N^{j+1}@>>>0\,.
\endalignedat
$$
By assumption, the ring $R_\p$ has global dimension $\le j$ whenever
$\codim\p\le c$, hence Lemma 11 shows that $\lau^j(N^{j+1})=0$.  The additivity
of $\f_R^j(M,-)$ now follows by applying $\lau^j$ to the long exact sequence.
The additivity of $\f_R^j(-,N)$ yields to a similar approach.

(b) For each $i$ the $R$--module $\Ext_R^i(M,N)$ is annihilated by
$(\ann M +\ann N)$, hence we have $\dim R-\dim \Ext_R^i(M,N)\ge
\dim R-\dim R/(\ann M +\ann N)>j$, and thus $\f_R^j(M,N)$ vanishes
in the indicated range by Lemma 10(b).
\QED
\enddemo

\demo{Proof of Proposition \rom{7}}
Lemma 10(b) shows that $\lau^h(R/\p)$ and $\lau^0(R)$ are positive, and that
$$
\left[\phi_R(R/\p,R/\p)(t)\right]_1
=(-1)^h\cdot\frac{\lau^h(R/\p)}{\lau^0(R)}\cdot\frac{\lau^h(R/\p)}{(1-t)^{d-2h}}
+O\left(\frac1{(1-t)^{d-2h-1}}\right)\,.
$$
With $\eta_h(\p)=\sum_{i=0}^h(-1)^i
\rank_{k(\p)}\!\left(\Ext^i_{R_\p}\!\left(k(\p),k(\p)\right)\right)$,
Lemmas 12(b) and 10(c) give
$$
\sum_{j=0}^h\frac{\f_R^j(R/\p,R/\p)}{(1-t)^{d-j}}=
\eta_h(\p)\cdot\frac{\lau^h(R/\p)} {(1-t)^{d-h}}\,.
$$

Thus, $\f_R^*(R/\p,R/\p)$ and $\phi^*_R(R/\p,R/\p)$ agree up to level $0$
if and only if $\lau^0(R)=\lau^0(R/\p)$.  In view of Lemma 10(c),
this is equivalent to the statement that $\p$ is the only prime of
$R$ with $\codim\p=0$ and that $R_\p$ is a field.

If $h>0$, then $\f_R^*(R/\p,R/\p)$ and $\phi^*_R(R/\p,R/\p)$ agree up to level
$h$ if and only if $\eta_h(\p)=0$. When the local ring $R_\p$ is regular the
rank of $\Ext^i_{R_\p}\!\left(k(\p),k(\p)\right)$ over $k(\p)$ is equal to
$\binom hi$, hence $\eta_h(\p)=0$. When $R_\p$ is not regular its embedding
dimension is at least $2$ because $\dim R_\p=h\ge1$, and then $\eta_h(\p)\ne0$
by Okiyama \cite{9\rm; Proposition 8}.
\QED
\enddemo

\proclaim{Lemma 13}
\rom{(a)} If $M$ and $N$ are finite graded $R$--modules, then for $a,b\in\BZ$
the sequences $\f_R^*(M(a),N(b))$ and $\phi^*_R(M(a),N(b))$ agree up to level
$n$ if and only if $\f_R^*(M,N)$ and $\phi^*_R(M,N)$ agree up to the same level.

\rom{(b)} If $\p$ is a homogeneous prime ideal such that $\codim\p=h$ and the
local ring $R_\p$ is Gorenstein, then $\f_R^*(R/\p,R)$ and $\phi^*_R(R/\p,R)$
agree up to level $h$.
\endproclaim

\demo{Proof}
(a) The canonical isomorphisms $\Ext^i_R(M(a),N(b))\cong\Ext^i_R(M,N)(b-a)$
show that
$$
\sum_{j=0}^n\frac{\f_R^j(M(a),N(b))}{(1-t)^{d-j}}=
t^{a-b}\cdot\sum_{j=0}^n\frac{\f_R^j(M,N)}{(1-t)^{d-j}}
$$
holds for each $n\ge 0$, while it is clear that
$\phi_R(M(a),N(b))(t)=t^{a-b}\cdot\phi_R(M,N)(t)$.

(b) By a direct computation as in the proof of Proposition 7, we get
$$
\left[\phi_R(R/\p,R)(t)\right]_1 =(-1)^h\cdot\frac{\lau^h(R/\p)}{(1-t)^{d-h}}
+O\left(\frac1{(1-t)^{d-h-1}}\right)\,.
$$
On the other hand, the Gorenstein hypothesis yields $\mu^i_R(\p,R)=\delta_{ih}$,
hence
\TagsOnRight
$$
\sum_{j=0}^h\frac{\f_R^j(R/\p,R)}{(1-t)^{d-j}}=
(-1)^h\cdot\frac{\lau^h(R/\p)}{(1-t)^{d-h}}\,.\tag"\qed"
$$
\TagsOnLeft
\enddemo
\pagebreak

\demo{Proof of Theorem \rom{8}}
Due to Lemmas 10(a) and 12(a), both sides of (4.$n$) are biadditive functions
on $\G(R)$ for $n=0,1,2$. Thus, it suffices to consider $M$ and $N$ which are
translates of quotients of $R$ by homogeneous prime ideals. By Lemma 13(a), we
can even assume that $M=R/\p$ and $N=R/\q$ with $\p,\q\in\Proj R$.  

If $\codim\p>n$, or $\codim\q>n$, or $\codim\p=n=\codim\q$ and $\p\ne\q$,
then the left hand side of (4.$n$) vanishes by Lemma 12(b), and its right
hand side vanishes by Lemma 10(b).  If $\codim\p=n$ and $\p=\q$ (respectively,
if $\p=0$ or if $\q=0$), then (4.$n$) holds by Proposition 7 (respectively, by
Proposition 6 or by Lemma 13(b)).

It remains to prove (4.2) when $\codim\s=1$, where $\s=\p$ or $\s=\q$.
As $R$ is factorial, $\s=(x)$ for a homogeneous $x\in R$,
hence by additivity and Lemma 13(a) we get the result from the
exact sequence $0@>>>R(-\deg x)@>>>R@>>>R/\s@>>>0$ and already settled cases.
\QED
\enddemo

\example{Remark 14}
Let $M$ be a finite graded module over a graded integral domain $R$.

If $R$ is Gorenstein in codimension $1$, then
$$
\lau^1(\Hom_R(M,R))-\lau^1(\Ext^1_R(M,R))=2\lau^1(R)\rank_R(M)-\lau^1(M)\,.
\tag6.1
$$
If $R$ is factorial, then in addition
$$
\gathered
\lau^2(\Hom_R(M,R))-\lau^2(\Ext^1_R(M,R))+\lau^2(\Ext^2_R(M,R))=\\
\biggl(1+2\frac{\lau^1(R)}{\lau^0(R)}\biggr)
\biggl(\lau^1(R)\rank_R(M)-\lau^1(M)\biggr)+\lau^2(M)\,.
\endgathered \tag6.2
$$

Indeed, the formulas above are obtained by rewriting (4.1) and (4.2) for
$N=R$ with the help of (5).  A factorial domain is Cohen-Macaulay in
codimension $2$ (due to normality), and thus is Gorenstein in this codimension
by a result of Murthy, cf\. \cite{4\rm; (3.3.19)}.  The argument for Lemma
12(a) shows that $\f_R^j(-,R)$ is an additive function on $\G(R)$ for $j=0,1,2$.
At this point we may retrace the preceding proof, ignoring any discussion
of primes $\q\ne0$.
\endexample

\head
4. Normal domains
\endhead

Our examples are built from the ring $R=K[X,Y,U,V]/(XV-YU)=K[x,y,u,v]$
with the standard grading, and the $R$--module $M=R/(u,v)\cong K[X,Y]$.

\proclaim{Example 15}
Equality \rom{(4.2)} may fail over a $3$-dimensional graded hypersurface
ring $R$ which is regular in codimension~$2$ (and hence is a normal domain).
\endproclaim

Over $K[X,Y,U,V]$ the quadratic form $XV-YU$ has a homogeneous {\it matrix
factorization\/}
$$
\pmatrix V & -Y\\ -U &X\endpmatrix\pmatrix X &Y\\U &V\endpmatrix=
\pmatrix XV-YU & 0\\0 & XV-YU\endpmatrix=
\pmatrix X &Y\\U &V\endpmatrix\pmatrix V & -Y\\ -U &X\endpmatrix\,,
$$
hence Eisenbud \cite{6\rm; (4.1)} provides the following minimal graded free
resolution of $M$:
$$
\dots@>>>R^2(-4)@>\pmatrix v &\!-y\\\!-u &x\endpmatrix>>\!
R^2(-3)@>\pmatrix x &y\\u &v\endpmatrix>>\!
R^2(-2)@>\pmatrix v &\!-y\\\!-u &x\endpmatrix>>\!
R^2(-1)@>\dsize{(u,v)}>>\!R@>>>0\,.
$$
\pagebreak

Applying to it the functor $\Hom_R(-,M)$ we get a complex of graded
$R$--modules
$$
0@>>>M@>\dsize{0}>>M^2(1)@>\pmatrix 0&0\\-y&x\endpmatrix>>
M^2(2)@>\pmatrix x&0\\y&0\endpmatrix>>
M^2(3)@>\pmatrix 0&0\\-y&x\endpmatrix>>
M^2(4)@>>>\dots\ .
$$
Direct computations based on the $M$--regularity of the sequence $x,y$ yield
$$
\gather
M\cong\Hom_R(M,M)\cong\Ext^1_R(M,M)\,;\\
\Ext^{2i}_R(M,M)\cong k(2i) \text{\quad for } i\ge1\,;\\
\Ext^{2i+1}_R(M,M)=0 \text{\quad for } i\ge1\,.
\endgather
$$
Accordingly, we get Laurent expansions
$$
\gather
H_R(t)=\frac{1-t^2}{(1-t)^4}=\frac2{(1-t)^3}-\frac{1}{(1-t)^2}\,;\\
H_M(t)=H_{\Hom_R(M,M)}(t)=H_{\Ext^1_R(M,M)}(t)=\frac1{(1-t)^2}\,;\\
H_{\Ext^{2}_R(M,M)}(t)=\frac1{t^2}=\sum_{j\ges0}(j+1)(1-t)^j\,.
\endgather
$$

It follows that for $M=N$ the expression on the left hand side of (4.2) is
equal to $0$, while the one on its right hand side is equal to $-1$.
\medskip

\proclaim{Example 16}
Equality \rom{(6.2)} may fail over a $4$-dimensional graded Cohen-Macaulay ring
$R'$ of type $2$ which is regular in codimension $3$ (and hence is a normal
domain).
\endproclaim

Let $R'=R[Z,W]/(xW-uZ,yW-vZ)=K[x,y,z,u,v,w]$ with the standard
grading, and set $M'=R'/(u,v,w)\cong K[X,Y,Z]$.  Since $w$ is $R'$--regular
and annihilates $M'$, for $S=R'/(w)\cong R[Z]/(uZ,vZ)$ and for each
$i$ we have $\Ext^i_{R'}(M',R')\cong\Ext^{i-1}_S(M',S)(1)$.  As
$M'=M[Z]$, and the minimal free resolution of $M[Z]$ over $R[Z]$ is ``the same
as'' that of $M$ over $R$, it is easy to get a beginning of a minimal
$S$--free resolution of $M'$ in the form
$$
S^4(-2)@>\pmatrix v & -y  & z & 0 \\ -u & x & 0 & z \endpmatrix>>
S^2(-1)@>\dsize{(u,v)}>>S@>>>0\,.
$$
Thus, we are interested in the homology of the complex
$$
0@>>>S@>{\tsize \alpha}=\pmatrix u\\v \endpmatrix>>S^2(1)
@>{\tsize \beta}=\pmatrix v&-u\\-y&x\\z&0\\0&z \endpmatrix>> S^4(2)\,.
$$

In degree $0$ we have $\Ker\alpha=(0\,:_S(u,v))=zS\cong M'(-1)$, hence
$\Ext^1_{R'}(M',R')\cong M'$.

In degree $1$ the complex is exact, hence $\Ext^2_{R'}(M',R')=0$.  Indeed,
writing $b\in S^2(1)$ in the form $(r_1+s_1z,r_2+s_2z)$ with uniquely
determined $r_j\in R$, we see that $\beta(b)=0$ implies
$$
\xalignat4
vr_1-ur_2&=0 &-yr_1+xr_2&=0 &z^2s_1&=0 &z^2s_2&=0\,.
\endxalignat
$$
 From the last two equations we get
$$
s_j\in(0\,:_S z^2)=((0\,:_S z)\,:_S z)=((u,v)\,:_S z)=(u,v)S\,,
$$
hence $b=(r_1,r_2)$. In view of the preceding example, the first two equations
yield an $a\in R$ such that $r_1=au$ and $r_2=av$, hence $b=a(u,v)\in\Im\alpha$.

We can now collect the relevant data on Laurent expansions:
$$
\gather
H_{R'}(t)=\frac{1-3t^2+2t^3}{(1-t)^6}=\frac3{(1-t)^4} - \frac2{(1-t)^3}\,;\\
H_{M'}(t)=H_{\Ext^1_{R'}(M',R')}(t)=\frac1{(1-t)^3}\,;\\
H_{\Hom_{R'}(M',R')}(t)=H_{\Ext^2_{R'}(M',R')}(t)=0\,.
\endgather
$$

They show that the left hand side of (6.2) is equal to $0$, and its right hand
side to $\frac13$.

\enddocument